\documentclass{article}
\usepackage{graphicx} % Required for inserting images
\usepackage{amsmath}
\usepackage{amsthm}
\usepackage{amsfonts} 
\usepackage{dsfont}
\usepackage{chngcntr}
\usepackage{natbib}
\usepackage{float}
\usepackage{lipsum}
\usepackage{multirow}
\usepackage{authblk}
\usepackage{authblk}
\usepackage{authblk}
\newtheorem{Proposition}{Proposition}
\begin{document}
\title{A Frequentist Approach to Change Point Detection: Methods and Applications}
\date{}

\author{Debanjana Datta}
%\footnotemark

\affil{Indian Statistical Institute, Bangalore}
\renewcommand\thefootnote{}\footnotetext{*Corresponding author: dbnjanad@gmail.com}

\maketitle
\begin{abstract}

In this paper we study the problem of change point detection in functional time series where the observations are allowed to vary on both sparse and dense support. We address the problem of  mean shift as well as the process volatility. Our methodology is based on the maximization of the conditional probability of change point given all the other parameters. Further, it has been proved that the proposed estimator is consistent.
\end{abstract}


\section{Introduction}
In recent times, the importance of functional data analysis lies in its ability to handle and extract meaningful information from complex phenomena occurring in biomedicine, finance, meteorology, environmental science, and others.Most of the works rely on the assumption of independent and identically distributed data, which is often not the case in observational data.\\ With the lapse of time, the inherent structure of the data may change due to the underlying stochastic process generating it.Hence, such an instance is of prior importance.

\section{Literature Review}
In functional data, our observations are random continuous trajectories generated from some   process and hence these are ultra-high-dimensional objects. In practice, these functions are observed at some finitely many grids in a compact domain.  The most common approach used to model functional observations is to project them onto some finite-dimensional space using some dimension reduction technique, namely Functional Principal Component Analysis (FPCA).\cite{ramsay2005functional} have enhanced the FDA
literature with  detailed discussions on several techniques and usefulness of FPCA.The functional autoregressive (FAR) model, developed by \cite{bosq2000linear},
plays a central role in modeling and predicting functional time series.\cite{gabrys2007portmanteau} proposed a test for independence and
identical distribution of functional observations.\\
For functional data, analysis has mainly been carried out under the assumption
of a single change point. Noted works in this area include 
\cite{horvath2010testing} where a test for equality of the functional autoregressive operator (FAR),
 \cite{aue2009break} where a test has been formulated to access the stability of volatilities and crossvolatilities of a process whereas \cite{berkes2009detecting} and 
tested the assumption of constant mean in a sequence of functional observations and   \cite{aue2009estimation} have dealt with the position of such a change.Further \cite{banerjee2018more} have provided a different test  for detecting the change point in mean in a sequence of independent as well as weakly dependent functional data.\\ 
\cite{aue2020testing} improvised a test for stationarity functional observations.Further, \cite{van2021nonparametric} proposed a non-parametric test of stationarity based on an explicit representation of the $L^{2}$
distance between the spectral density operator of a nonstationary process and
its best $L^{2}$ approximation by a spectral density operator corresponding to a
stationary process.
A different version of functional stationarity
tests, based on time domain methodology involving cumulative sum statistics  was given in \cite{horvath2014testing}.
Following the Bayesian perspective pioneering works have been done by \cite{chernoff1964estimating},\cite{smith1975bayesian} and \cite{carlin1992hierarchical} and recently by \cite{li2021bayesian}. \\
The rest of the paper is organized as follows. Section \ref{3} illustrates the FAR model  and its representation in a dynamic linear model (DLM) form, which is further used for inferential purposes. Algorithms have been proposed in section \ref{4} whereas in Section \ref{5}  simulation results have been provided under different setups. Data analysis have been carried out in Section \ref{6}.

\section{Model}{\label{3}}
Let  $Y_{1}$,\ldots,$Y_{T}$ be a time ordered sequence of random functions in \ $L^{2}(\mathcal{T})$ where $\mathcal{T}$ is a compact indexing set. For brevity, consider  $\mathcal{T}$ =[0,1].\\A Functional Autoregressive process of order 1 (FAR(1))  is given by
\begin{math}
    Y_{t}=\Psi (Y_{t-1})+\epsilon_{t}\
\end{math}
;\ t = 1,\ldots, T
where
\begin{math}
\Psi : L^{2} \rightarrow L^{2}
\end{math} is an autoregressive operator and 
\begin{math}
\{\epsilon_{t}\}_{t\leq T}
\end{math}
is a sequence of independently and identically distributed(iid) random errors in \ $L^{2}(\mathcal{T})$.
We restrict our attention to integral operators only defined by \begin{math}
    \Psi(Y_{t})(k)= \int_{0}^{1} \psi (k,s) Y_{t-1}(s)\,ds.
\end{math}
These functional observations $Y_{t}$ are not observed directly, but rather via  discrete samples of each curve, and typically with measurement error.\\ 
Suppose we observe 
\begin{math}
y_{i,t} \in \mathcal{R}
\end{math} 
\ sampled with noise
$v_{i,t}$ from 
\begin{math}
Y_{t} \in L^{2}:
\end{math}\\
\begin{math}
    y_{i,t}=Y_{t}(k_{i,t})+ v_{i,t}
\end{math}
; t = 1,\ldots,T;
where 
\begin{math}
  {\tau_t} =\{ k_{1,t},\ldots, k_{m_{t},t}\}
\end{math}
are observation points of $Y_{t}$ and $v_{i,t}$ is a mean zero measurement error with finite variance.\\

\begin{Proposition}({\cite{kowal2019functional}})
 Let \begin{math}
     X_{t}-\mu=\sum_{i=1}^{p} \Psi_i(X_{t-i}-\mu)+\epsilon_{t}, 
 \end{math} and suppose that we observe
 \begin{math}
     x_{t}=X_{t}+\nu_{t},
 \end{math} where $\{\epsilon_t\}$ and $\{\nu_t\}$ are two independent white noise processes. Then, the observable process $\{x_t\}$ follows a Functional Autoregressive Moving Average process of order(p,p) (FARMA(p,p)).
 \end{Proposition}

Let $\alpha_{t}(.)$ = $Y_{t}(.)$-$\mu(.)$. Hence the two level hierarchical
model is given as:\\

\begin{equation*}
\begin{split}
 y_{i,t} & =\mu(k_{i,t})+\alpha_{t}(k_{i,t})+  v_{i,t} \\
 \alpha_{t}(k) & =\int_{0}^{1} \psi (k,s) \alpha_{t-1}(s)\,ds + \epsilon_{t}(\tau) \ , \forall \ k \in \mathcal{T}
\end{split}
\end{equation*}
\\
where we assume \begin{math}
\{v_{i,t}\} \end{math} and \begin{math} \{\epsilon_{t}\} \end{math} are mutually independent sequences.
\\

This two level hierarchical model bridges that gap between model misspecification and hence increases forecast accuracy.
For implementation of model(1), a finite set of evaluation points \begin{math} \tau_{e} = \{k_{1},\ldots,k_{M}\} \in \mathcal{T} \end{math} is to be selected. Assume, \begin{math}
    \tau_{t} \subseteq \tau_{e} \forall
\end{math}\ t. \\
Using quadrature methods approximate the integral in (1) accurately by choosing M large and $\tau_{e}$ to be dense in \begin{math}\mathcal{T} : \end{math}
\begin{math}
     \int_{0}^{1} \psi(k,s)\alpha_{t-1}(s)\,ds  \approx  (\psi(k,k_{1}),\ldots,\psi(k,k_{{M}})) \mathbf{Q}\mathbf{\alpha_{t-1}}
\end{math} 
\\
where $\mathbf{Q}$ is a known quadrature weight matrix and 
\begin{math} 
\mathbf{\alpha_{t-1}} = (\alpha_{t-1}(k_{1}),\ldots,\alpha_{t-1}(k_{M})) \end{math}
\\ 
Let $\mathbf{Z_t}$ is the $m_{t}\times M$ incidence matrix that identifies the observation points at time t. We can have a state space formulation of model(1) as a dynamic linear model (\cite{west2006bayesian}) in $\mathbf{\alpha_{t}(.)}$:

\begin{align*}
    \mathbf{y_{t}}   &   = \mathbf{Z_t}\mathbf{\mu} + \mathbf{Z_t}\mathbf{\alpha_{t}} + \mathbf{\nu_{t}} ; & \mathbf{\nu_{t}}\stackrel{indep} \sim N(0,\sigma_\nu^{2}\mathds{I}_{m_t})  &   ,&t=1,2,..,T \\
    \mathbf{\alpha_{t}} & =\mathbf{\Psi} \mathbf{Q}\mathbf{\alpha_{t-1}}+\mathbf{\epsilon_{t}}\vspace{5mm}\ ;&\mathbf{\epsilon_{t}}\stackrel{indep} \sim N(0,\mathbf{K_\epsilon)}& ,& t=2,....,T\\
    \mathbf{\alpha_{1}}& \sim N(0,K_{\epsilon})
\end{align*}\\
where 
\begin{math}\mathbf{y_{t}} = (y_{1,t},\ldots,y_{m_{t},t}) \end{math} ,
\begin{math}\mathbf{\mu} = (\mu({k_1}),\ldots,\mu({k_M}))\end{math}, 
\begin{math}
\mathbf{\Psi}= \{\psi(k_{i},k_{j})\}_{i,j=1}^{M}
\end{math} and 
\begin{math}
    \mathbf{K_{\epsilon}}= \{K_{\epsilon}(k_{i},k_{j})\}_{i,j=1}^{M}
\end{math}\\

The change point   \begin{math} \tau \end{math} in the mean function is defined as follows:

\[
\mathbf{y_t} = \begin{cases}
    \mathbf{Z_t}\mathbf{\mu_1} + \mathbf{Z_t}\mathbf{\alpha_t} + \mathbf{\nu_t} ; & t=1,2,..,\tau \\
    \mathbf{Z_t}\mathbf{\mu_2} + \mathbf{Z_t}\mathbf{\alpha_t} + \mathbf{\nu_t} ; & t=\tau+1,..,T
\end{cases}
\quad \text{; } \mathbf{\nu_t} \stackrel{indep}{\sim} N(0,\sigma_\nu^{2}\mathds{I}_{m_t})
\]

Similarly, in case of process volatility,  the change point   \begin{math} \tau \end{math} is given by

\[
\mathbf{y_t} = \begin{cases}
    \mathbf{Z_t}\mathbf{\mu} + \mathbf{Z_t}\mathbf{\alpha_t} + \mathbf{\nu_t} ; & \mathbf{\nu_t} \stackrel{indep}{\sim} N(0,\sigma_1^{2}\mathds{I}_{m_t}) \quad \text{for } t=1,2,..,\tau \\
    \mathbf{Z_t}\mathbf{\mu} + \mathbf{Z_t}\mathbf{\alpha_t} + \mathbf{\nu_t} ; & \mathbf{\nu_t} \stackrel{indep}{\sim} N(0,\sigma_2^{2}\mathds{I}_{m_t}) \quad \text{for } t=\tau+1,..,T
\end{cases}
\]

\section{Algorithms}{\label{4}}
\subsection{Mean Problem}
{
\begin{itemize}
    \item For each
    $\tau \in \{2, \dots, T-1\}$ 
    \begin{itemize}
    \item Estimate $\mu_1$ and $\mu_2$ by the method of least squares.
    \item Fit a spline $\{y_t -\hat{\mu}\}_{t=1}^T $  to get the $\{\hat{\alpha}_t\} _{t=1}^{T}$
    \item Obtain the MLE estimate $\hat{\sigma}_\nu^{2}$ .
    \item Estimate $\Psi$ as a function of $\{\hat{\alpha}_t\} _{t=1}^{T}$ .
    
 $\Psi$ is based on the parametrization 
\begin{math}
\Psi=B_{\Psi}\Theta_{\Psi}B_{\Psi}^{'}
\end{math} where $B_{\Psi}$ is a matrix of cubic B-spline basis functions and
\begin{math}
\theta_{\Psi} = vec(\Theta_{\Psi}) 
\end{math} is sampled from 
\begin{math}
\mathcal{N}(A_{\psi}a_{\psi},A_{\psi}),
\end{math}
where 
\[
A_{\psi}^{-1}=\Omega_{\psi}+[(B_{\psi}^{'}Q)\{\sum_{t=2}^{T}\hat{\alpha}_{t-1}\hat{\alpha}_{t-1}^{'}\}(B_{\psi}^{'}Q)^{'}]\otimes[(B_{\psi}^{'}B_{\psi})] \]
\[a_{\psi}=vec(B_{\psi}^{'}\sum_{t=2}^{T}\hat{\alpha_{t}}\hat{\alpha}_{t-1}^{'}(B_{\psi}^{'}Q)^{'})
\]

    \item Re-estimate the $\{{\alpha}_t\} _{t=1}^{T}$ and hence ${\mu}_1$, ${\mu}_2$, and ${\sigma}_\nu^{2}$ according to the dynamic linear model structure.
    \end{itemize}
    \item Find that particular value of $\tau$ which maximizes the conditional probability [$\tau|$\ldots].
\end{itemize}

\subsection{Volatility Problem}

\begin{itemize}
        \item Estimate $\mu$ by the method of least squares.
        \item Estimate $\{\alpha_t\}_{t=1}^{T}$ by fitting a spline to $\{y_t -\hat{\mu}\}_{t=1}^T $.
        \item For each  \begin{math}
        \tau \in \{2,..,T-1\} \end{math} 
        \begin{itemize}
            \item 
       Estimate $\sigma_{1}^{2}$ and $\sigma_{2}^{2}$  by the method of maximum likelihood.
         \item Estimate $\Psi$ as a function of $\{\hat{\alpha}_t\} _{t=1}^{T}$ similarly as the previous one.
    
    \item Re-estimate $\{{\alpha}_t\} _{t=1}^{T}$  according to the dynamic linear model structure. 
        \item Hence, re-estimate $\mu$, $\sigma_{1}^{2}$ and $\sigma_{2}^{2}$.
         \end{itemize}
        \item Find that particular value of $\tau$ which maximizes the conditional probability [$\tau|$\ldots].
    \end{itemize}
  
\section{Simulation Studies}{\label{5}}

For the mean problem, we have considered  ($\mu_1(u),\mu_2(u)$) = ($\frac{1}{10} u^3 \sin(2\pi u)$, $\frac{1}{10} u^3 \sin(2\pi u)+u^{2}$) and ($sin(u),cos(u)$). The measurement errors  $\nu_{i,t} \overset{\text{iid}}{\sim} N(0, \sigma^2_\nu)$ with $\sigma_\nu = 0.001$ and the FAR(1) kernel used is Bimodal-Gaussian kernel, $\psi(s, u) \propto \frac{0.75}{\pi(0.3)(0.4)} \exp\left\{ -\frac{(s - 0.2)^2}{(0.3)^2} - \frac{(u - 0.3)^2}{(0.4)^2}\right\} + \frac{0.45}{\pi(0.3)(0.4)} \exp\left\{-\frac{(s - 0.7)^2}{(0.3)^2} - \frac{(u - 0.8)^2}{(0.4)^2}\right\}$; rescaled according to a pre-specified squared norm, $C_\psi = \int \int \psi^2(s, u) \,ds \,du$, with $C_\psi < 1$ for stationarity. We select $C_{\psi} = 0.4$ for the simulation. We vary the sample size from small T = 50 to large T = 200
and the grid sizes from sparse m = 30 to dense m = 100. The estimated change points exactly matches with the true ones for all sample sizes.\\
For the volatility problem, we have considered the mean function  $\mu(u) = \frac{1}{10} u^3 \sin(2\pi u)$ and the same kernel as in the previous one.  It is noted that for every T fixed increasing grid size gives accurate results.

\begin{table}[H]
\centering
\caption{$\sigma_1 = 0.001$ , $\sigma_2 = 0.002$}
  \begin{tabular}{|c|c||c|c|c|}
    \hline
    T & m & CP at $\left\lceil \frac{T}{4} \right\rceil$ & CP at T/2 & CP at  $\left\lceil \frac{3T}{4} \right\rceil$ \\
    \hline
    \multirow{3}{*}{50} 
      & 30 & 3  & 17 & 39 \\
      & 50 & 45 & 25 & 36 \\
      & 100 & 14 & 25 & 37 \\
    \hline
    \multirow{3}{*}{100} 
      & 30 & 25 & 49 & 74 \\
      & 50 & 25 & 49 & 76 \\
      & 100 & 26 & 50 & 75 \\
    \hline
    \multirow{3}{*}{200} 
      & 30 & 48 & 109 & 147 \\
      & 50 & 51 & 96  & 150 \\
      & 100 & 50 & 100 & 147 \\
    \hline
  \end{tabular}
\end{table}

\begin{table}[H]
\centering
\caption{$\sigma_1 = 0.001$ , $\sigma_2 = 0.01$}
  \begin{tabular}{|c|c||c|c|c|}
    \hline
    T & m & CP at $\left\lceil \frac{T}{4} \right\rceil$ & CP at T/2 & CP at $\left\lceil \frac{3T}{4} \right\rceil$ \\
    \hline
    \multirow{3}{*}{50} 
      & 30 & 13  & 25 & 38 \\
      & 50 & 13 & 25 & 38 \\
      & 100 & 13 & 25 & 38 \\
    \hline
    \multirow{3}{*}{100} 
      & 30 & 25 & 50 & 71 \\
      & 50 & 25 & 50 & 75 \\
      & 100 & 25 & 50 & 75 \\
    \hline
    \multirow{3}{*}{200} 
      & 30 & 49 & 99 & 148 \\
      & 50 & 50 & 100  & 150 \\
      & 100 & 50 & 100 & 150 \\
    \hline
  \end{tabular}
\end{table}

\section{Real Data Analysis}{\label{6}}
For the first real data analysis we compare the performance of the proposed methodology with  \cite{banerjee2018more} and \cite{berkes2009detecting}. We have studied the average daily temperatures of central England for 228 years, from 1780 to 2007. This data can be viewed as 228 curves with 365 measurements on each curve. We shall study this dataset under mean shift problem.
\begin{table}[H]
\centering
\caption{Performance comparison of the algorithm for Central England temperature data}
\label{tab:performance_comparison}
\begin{tabular}{|c||c|c|c|}
\hline
\textbf{Year} & \textbf{Estimated CP\footnotemark} & \textbf{Estimated CP\footnotemark} & \textbf{Estimated CP} \\
\hline
1780-2007 & 1926 & 1926 & 1919 \\
\hline
\end{tabular}
\end{table}
\footnotetext[1]{\citeauthor{berkes2009detecting}}
\footnotetext[2]{\citeauthor{banerjee2018more}}

\begin{figure}[h]
\centering
\includegraphics[width=0.45\textwidth]{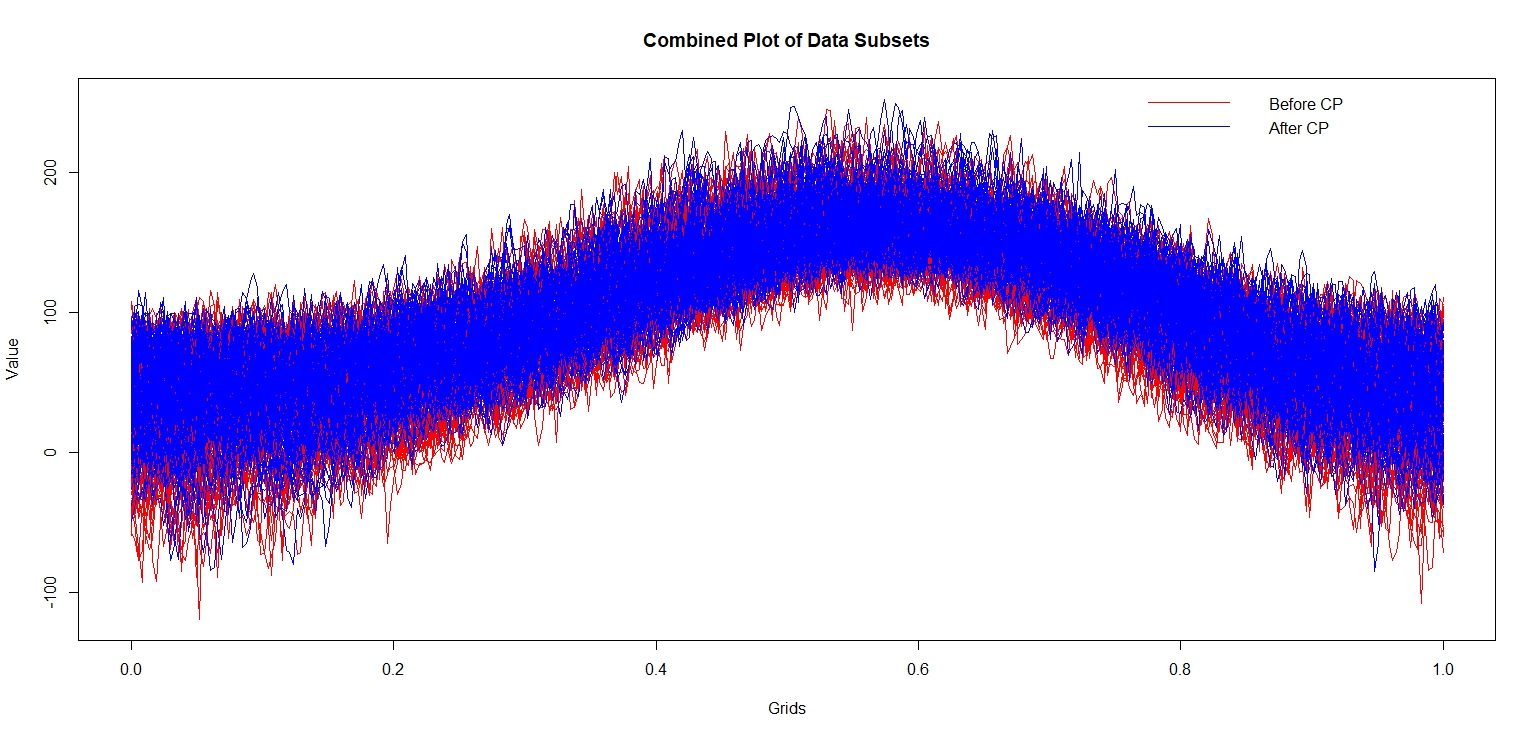}
\end{figure}
}

For the second application, logarithm of the closing prices of AdaniPower stock from 7th Jan'19 to 1st Jan'21 is considered on a weekly basis at a definite interval of time.
Our estimated changepoint for the AdaniPower stock is the week of 20th July 2020.

\renewcommand\thefootnote{}\footnote{CP: Change Point}
\addtocounter{footnote}{-1}

\bibliography{References}
\section*{Appendix: Proof}

Model:   
\begin{align*}
    \mathbf{y_{t}}   &   = \mathbf{Z_t}\mathbf{\mu_1} + \mathbf{Z_t}\mathbf{\alpha_{t}} + \mathbf{\nu_{t}} ; & \mathbf{\nu_{t}}\stackrel{indep} \sim N(0,\sigma_\nu^{2}\mathds{I}_{m_t})  &   ,&t=1,2,..,\tau \\
        \mathbf{y_{t}}   &   = \mathbf{Z_t}\mathbf{\mu_2} + \mathbf{Z_t}\mathbf{\alpha_{t}} + \mathbf{\nu_{t}} ; & \mathbf{\nu_{t}}\stackrel{indep} \sim N(0,\sigma_\nu^{2}\mathds{I}_{m_t})  &   ,&t=\tau+1,..,T\\
\end{align*} 

Let, $\hat{\mu_1}$ and $\hat{\mu_1}^{E}$ be the estimates of $\mu_1$ under the true change point $\tau$ and the estimated change point $\hat{\tau}$ respectively.\\
Similarly,  $\hat{\mu_2}$ and $\hat{\mu_2}^{E}$ be the estimates of $\mu_2$ under the true change point $\tau$ and the estimated change point $\hat{\tau}$ respectively.
\begin{enumerate}
    \item Estimating $\mu_1$ and $\mu_2$ as a smooth mean of $\{\mathbf{y_t}\}_{t=1}^{T}$ under $\tau$ under the parametrization \begin{math}
        \mu_1=B\theta_1 ;
        \mu_2=B\theta_2.
    \end{math}
\begin{align*}
\left.
\begin{array}{ll}
    \mathbf{y_{t}}   = \mathbf{Z_t}\mathbf{\mu_1} + \mathbf{\nu_{t}}, & \mathbf{\nu_{t}}\stackrel{indep} \sim N(0,\sigma_\nu^{2}\mathds{I}_{m_t}), \quad t=1,2,\ldots,\tau \\
    \mathbf{y_{t}}   = \mathbf{Z_t}\mathbf{\mu_2} + \mathbf{\nu_{t}}, & \mathbf{\nu_{t}}\stackrel{indep} \sim N(0,\sigma_\nu^{2}\mathds{I}_{m_t}), \quad t=\tau+1,\ldots,T
\end{array}
\right\}\\
\quad \Longleftrightarrow \quad 
\left\{
\begin{array}{ll}
    \mathbf{y_{t}}   = \mathbf{Z_t}\mathbf{B\theta_1 } + \mathbf{\nu_{t}}, & \mathbf{\nu_{t}}\stackrel{indep} \sim N(0,\sigma_\nu^{2}\mathds{I}_{m_t}), \quad t=1,2,\ldots,\tau \\
    \mathbf{y_{t}}   = \mathbf{Z_t}\mathbf{B\theta_2 } + \mathbf{\nu_{t}}, & \mathbf{\nu_{t}}\stackrel{indep} \sim N(0,\sigma_\nu^{2}\mathds{I}_{m_t}), \quad t=\tau+1,\ldots,T
\end{array}
\right.
\end{align*}

By the method of least squares, under $\tau$:
\begin{align*}
    \hat{\theta}_{1} &= (A^{T}_{1}A_1)^{-1}A^{T}_1Y_{1}   \\ 
    \hat{\theta}_{2} &= (A^{T}_{2}A_2)^{-1}A^{T}_{2}Y_{2}
\end{align*}
where
\[
\begin{array}{cc}
Y_1 =
\begin{bmatrix}
\mathbf{y_1} \\
\mathbf{y_2} \\
\vdots \\
\mathbf{y_\tau}
\end{bmatrix}
;
Y_2 =
\begin{bmatrix}
\mathbf{y_{\tau+1}} \\
\mathbf{y_{\tau+2}} \\
\vdots \\
\mathbf{y_T}
\end{bmatrix}
\end{array}
\]
and 
\[
\begin{array}{cc}
A_1 =
\begin{bmatrix}
\mathbf{Z_1B} \\
\mathbf{Z_2B} \\
\vdots \\
\mathbf{Z_\tau B}
\end{bmatrix}
;
A_2 =
\begin{bmatrix}
\mathbf{Z_{\tau+1}B} \\
\mathbf{Z_{\tau+2}B} \\
\vdots \\
\mathbf{Z_T B}
\end{bmatrix}
\end{array}
\]
Also, $$\hat{\theta}_{1} \xrightarrow{P} \theta_{1}$$  and 
    $$\hat{\theta}_{2} \xrightarrow{P} \theta_{2}$$.

By the method of least squares, under $\hat{\tau}$\\
\begin{math}
    \hat{\theta}_{1}^E=(A_1^{*\,T}A_1^{*})^{-1}A^{*\,T}_1Y_{1}^{*} \\
    \hat{\theta}_{2}^E=(A_2^{*\,T}A_2^{*})^{-1}A^{*\,T}_2Y_{2}^{*}\\
\end{math}
where
\[
\begin{array}{cc}
Y_1^{*} =
\begin{bmatrix}
\mathbf{y_1} \\
\mathbf{y_2} \\
\vdots \\
\mathbf{y_{\hat{\tau}}}
\end{bmatrix}
;
Y_2^{*} =
\begin{bmatrix}
\mathbf{y_{\hat{\tau}+1}} \\
\mathbf{y_{\hat{\tau}+2}} \\
\vdots \\
\mathbf{y_T}
\end{bmatrix}
\end{array}
\]
and 
\[
\begin{array}{cc}
A^{*}_1 =
\begin{bmatrix}
\mathbf{Z_1B} \\
\mathbf{Z_2B} \\
\vdots \\
\mathbf{Z_{\hat{\tau}} B}
\end{bmatrix}
;
A^{*}_2 =
\begin{bmatrix}
\mathbf{Z_{\hat{\tau}+1}B} \\
\mathbf{Z_{\hat{\tau}+2}B} \\
\vdots \\
\mathbf{Z_T B}
\end{bmatrix}
\end{array}
\]

\end{enumerate}

Consider these two quantities,
\begin{math}
    ||B\hat{\theta}_{1} - B\hat{\theta}_{1}^{E} ||_{2}^{2} \; \text{and} \;||B\hat{\theta}_{2} - B\hat{\theta}_{2}^{E} ||_{2}^{2}
\end{math}.

We know, \begin{math}
    \hat{\theta}_{1} \xrightarrow{}\theta_1 ;
    \hat{\theta}_{1}^{E} \xrightarrow{} \theta_1-(A_1^{*\,T}A_1^{*})^{-1}A^{*\,T}_1Z(\theta_2-\theta_1)
\end{math} \\
where 
\[
Z=
\begin{bmatrix}
\mathbf{0} \\
\vdots \\
\mathbf{0} \\
\mathbf{Z_{\tau+1} B} \\
\vdots \\
\mathbf{Z_{\hat{\tau}} B}
\end{bmatrix}
\]

Therefore, 
\begin{math}
    ||B\hat{\theta}_{1} - B\hat{\theta}_{1}^{E} ||_{2}^{2}=||(A_1^{*\,T}A_1^{*})^{-1}A^{*\,T}_1Z(\theta_2-\theta_1)||_{2}^{2}
\end{math}\\
Consider, 
\begin{math}
    Z_{t}^{T}Z_{t} =I_{M} \; \;  \forall \; t\\
\end{math}
Then,
\begin{math}
    ||(A_1^{*\,T}A_1^{*})^{-1}A^{*\,T}_1Z(\theta_2-\theta_1)||_{2}^{2}=||(\frac{\hat{\tau}-\tau}{\hat{\tau}})(\mu_2-\mu_1)||_{2}^{2} =(\frac{\hat{\tau}-\tau}{\hat{\tau}})^{2}||\mu_2-\mu_1||_{2}^{2}
\end{math}
Hence, 
\begin{math}
    ||B\hat{\theta}_{1} - B\hat{\theta}_{1}^{E} ||_{2}^{2} \xrightarrow{} (\frac{\hat{\tau}-\tau}{\hat{\tau}})^{2}||\mu_2-\mu_1||_{2}^{2}
\end{math}

We know, \begin{math}
    \hat{\theta}_{2} \xrightarrow{}\theta_2 ;
    \hat{\theta}_{2}^{E} \xrightarrow{} \theta_2-(A_2^{*\,T}A_2^{*})^{-1}A^{*\,T}_2Z^{*}(\theta_1-\theta_2)
\end{math} \\

where 
\[
Z^{*}=
\begin{bmatrix}
\mathbf{Z_{\hat{\tau}+1}B} \\
\vdots \\
\mathbf{Z_{{\tau}}B} \\
\mathbf{0} \\
\vdots \\
\mathbf{0}
\end{bmatrix}
\]

Therefore,
\begin{align*}
    ||B\hat{\theta}_{2} - B\hat{\theta}_{2}^{E} ||_{2}^{2} &= ||(A_2^{*\,T}A_2^{*})^{-1}A^{*\,T}_2Z^{*}(\theta_1-\theta_2)||_{2}^{2} \\
    &= ||(A_2^{*\,T}A_2^{*})^{-1}A^{*\,T}_2Z^{*}(\theta_2-\theta_1)||_{2}^{2}
\end{align*}
Consider,
\begin{align*}
    Z_{t}^{T}Z_{t} &= I_{M} \quad \forall \; t
\end{align*}
Then,
\begin{align*}
    ||(A_2^{*\,T}A_2^{*})^{-1}A^{*\,T}_2Z^{*}(\theta_2-\theta_1)||_{2}^{2} &= \left|\left|\left(\frac{\tau-\hat{\tau}}{T-\hat{\tau}}\right)(\mu_2-\mu_1)\right|\right|_{2}^{2} \\
    &= \left(\frac{\tau-\hat{\tau}}{T-\hat{\tau}}\right)^{2}||\mu_2-\mu_1||_{2}^{2}
\end{align*}
Hence,
\begin{align*}
    ||B\hat{\theta}_{2} - B\hat{\theta}_{2}^{E} ||_{2}^{2} &\xrightarrow{} \left(\frac{\tau-\hat{\tau}}{T-\hat{\tau}}\right)^{2}||\mu_2-\mu_1||_{2}^{2}
\end{align*}

Therefore we have proved,\\
\begin{align*}
    ||B\hat{\theta}_{1} - B\hat{\theta}_{1}^{E} ||_{2}^{2} &\xrightarrow{} \left(\frac{\hat{\tau}-\tau}{\hat{\tau}}\right)^{2}||\mu_2-\mu_1||_{2}^{2} \\
    ||B\hat{\theta}_{2} - B\hat{\theta}_{2}^{E} ||_{2}^{2} &\xrightarrow{} \left(\frac{\tau-\hat{\tau}}{T-\hat{\tau}}\right)^{2}||\mu_2-\mu_1||_{2}^{2}
\end{align*}

\[
\hat{\sigma_\nu^2} = \frac{\sum_{t=1}^{\tau} \sum_{i=1}^{m_t}(y_{i,t}-\widehat{\mu_1}(k_{i,t})-\widehat{\alpha}_t(k_{i,t}))^2 +   \sum_{t=\tau+1}^{T} \sum_{i=1}^{m_t}(y_{i,t}-\widehat{\mu_2}(k_{i,t})-\widehat{\alpha}_t(k_{i,t}))^2   }{\sum_{t=1}^T m_t}
\]

\[
\hat{\sigma_\nu^2}^{E} = \frac{\sum_{t=1}^{\hat{\tau}} \sum_{i=1}^{m_t}(y_{i,t}-\widehat{\mu_1}^{E}(k_{i,t})-\widehat{\alpha}^{E}_t(k_{i,t}))^2 +   \sum_{t=\hat{\tau}+1}^{T} \sum_{i=1}^{m_t}(y_{i,t}-\widehat{\mu_2}^{E}(k_{i,t})-\widehat{\alpha}^{E}_t(k_{i,t}))^2   }{\sum_{t=1}^T m_t}
\]

Hence,

\[
\hat{\sigma_\nu^2} = \frac{\sum_{t=1}^{\tau} ||y_t-Z_t \hat{\mu_1}-Z_t\hat{\alpha_t}||_{2}^{2} +   \sum_{t=\tau+1}^{T} ||y_t-Z_t \hat{\mu_2}-Z_t\hat{\alpha_t}||_{2}^{2}  }{\sum_{t=1}^T m_t}
\]

and
\[
\hat{\sigma_\nu^2}^{E} = \frac{\sum_{t=1}^{\tau} ||y_t-Z_t \hat{\mu_1}^{E}-Z_t\hat{\alpha_t}^{E}||_{2}^{2} +   \sum_{t=\tau+1}^{T} ||y_t-Z_t \hat{\mu_2}^{E}-Z_t\hat{\alpha_t}^{E}||_{2}^{2}  }{\sum_{t=1}^T m_t}
\]

Consider $\hat{\tau}=\tau$ :
\[
\hat{\sigma_\nu^2}^{E} = \frac{\sum_{t=1}^{{\tau}} ||y_t-Z_t \hat{\mu_1}-Z_t\hat{\alpha_t}||_{2}^{2} +   \sum_{t={\tau}+1}^{T} ||y_t-Z_t \hat{\mu_2}-Z_t\hat{\alpha_t}||_{2}^{2} + \sum_{t=1}^T ||Z_t\hat{\alpha_t}-Z_t\hat{\alpha_t}^{E} ||_{2}^{2}}{\sum_{t=1}^T m_t}
\]
\[
=\frac{\sum_{t=1}^{\tau} ||y_t-Z_t \hat{\mu_1}-Z_t\hat{\alpha_t}||_{2}^{2} +   \sum_{t=\tau+1}^{T} ||y_t-Z_t \hat{\mu_2}-Z_t\hat{\alpha_t}||_{2}^{2} }{\sum_{t=1}^T m_t}+\frac{\sum_{t=1}^T ||Z_t\hat{\alpha_t}-Z_t\hat{\alpha_t}^{E} ||_{2}^{2}}{\sum_{t=1}^T m_t}
\]
\[
=\hat{\sigma_\nu^2}+K
\]

Consider $\hat{\tau} < \tau $ :
\[
\begin{split}
\hat{\sigma_\nu^2}^{E} = \frac{\sum_{t=1}^{{\hat{\tau}}} ||y_t-Z_t \hat{\mu_1}-Z_t\hat{\alpha_t}||_{2}^{2} +   \sum_{t={\hat{\tau}}+1}^{\tau} ||y_t-Z_t \hat{\mu_2}-Z_t\hat{\alpha_t}||_{2}^{2} +   \sum_{\tau+1}^T ||y_t-Z_t\hat{\mu_2}-Z_t\hat{\alpha}_t||_{2}^{2} +}{\sum_{t=1}^T m_t}+ &\\
\frac{ \sum_{t=1}^T ||Z_t\hat{\alpha}_t-Z_t\hat{\alpha}^{E}_t||_{2}^{2}+\sum_{t=1}^{\hat{\tau}}||Z_t\hat{\mu_1}^{E}-Z_t\hat{\mu_1}||_{2}^{2}+ \sum_{t=\hat{\tau}+1}^{\tau}||Z_t\hat{\mu_1}-Z_t\hat{\mu_2}^{E}||_{2}^{2}+\sum_{{\tau}+1}^{T}||Z_t\hat{\mu_{2}}^{E}-Z_t\hat{\mu_2}||_{2}^{2}}{\sum_{t=1}^T m_t}
& \\ 
\hat{\sigma_\nu^2} +K+A_1 ; \quad  A_1=\frac{\sum_{t=1}^{\hat{\tau}}||Z_t\hat{\mu_1}^{E}-Z_t\hat{\mu_1}||_{2}^{2}+ \sum_{t=\hat{\tau}+1}^{\tau}||Z_t\hat{\mu_1}-Z_t\hat{\mu_2}^{E}||_{2}^{2}+\sum_{\tau+1}^{T}||Z_t\hat{\mu_{2}}^{E}-Z_t\hat{\mu_2}||_{2}^{2}}{\sum_{t=1}^T m_t}
\end{split}
\]

Consider $\hat{\tau} > \tau $ :
\[
\hat{\sigma_\nu^2}^{E} =\hat{\sigma_\nu^2} +K+A_2 ; \quad  A_2=\frac{\sum_{t=1}^{\tau}||Z_t\hat{\mu_1}-Z_t\hat{\mu_1}^{E}||_{2}^{2}+ \sum_{t={\tau+1}}^{\tau}||Z_t\hat{\mu_2}-Z_t\hat{\mu_1}^{E}||_{2}^{2}+\sum_{\hat{\tau}+1}^{T}||Z_t\hat{\mu_{2}}^{E}-Z_t\hat{\mu_2}||_{2}^{2}}{\sum_{t=1}^T m_t}
\]

\end{document}